\documentclass[12pt]{article}
\title{Factorization of Symmetric Matrices and Actions on $\Sigma_k$}
\author{Baohua FU}
\usepackage{amsmath,amssymb,amsthm}
\chardef\bslash=`\\
\newtheorem{Thm}{Theorem}

\newtheorem{Prop}[Thm]{Proposition}

\newtheorem{Rque}{Remark}

\begin{document}
\maketitle
\begin{abstract}
In this note, firstly we give an easy proof of the factorization of symmetric matrices (see \cite{Mos}), and then use it to prove the well-known 
fact that the automorphism group of a non-degenerate symmetric bilinear form $Q$ acts transitively on the locus of isotropic subspaces $\Sigma_k(Q)$.
\end{abstract}

Let $V$ be an $m$-dimensional complex vector space.
\begin{Prop}
Every symmetric matrix $Q \in End(V)$ is of the form $P^TP$, for some $P \in End(V)$.
\end{Prop}
\begin{proof}
We will not make any difference between symmetric matrices and symmetric bilinear forms on $V$. 
If $Q \neq 0$, then there exists some $u_1 \in V$ s.t. $Q(u_1,u_1) \neq 0$. Consider the hyperplane
$H_1 = \{u \in V | Q(u_1,u) = 0\}$. If $Q$ is not identically zero on $H_1$, then we can chose some $u_2 \in H_1$, s.t.
$Q(u_2, u_2) \neq 0$, and then $Q(u_2,u_1) =0$. Then we consider $H_2 =\{u \in H_1 | Q(u_2,u) = 0\}$, and we continue this procedure. At the end, 
we will get some vectors $v_1,\dots,v_l,v_{l+1},\dots, v_m$, such that 
$$Q(v_i,v_j) = \delta_{ij}, \, for\, 1 \leq i,j \leq l ; \quad Q(v_i,v_j) = 0, \, if\, max\{i,j\} \geq l+1. $$
In this basis, $Q$ is of the form $\,diag(I_l,0)$, i.e. there exists some invertible matrix $A$ such that $A^TQA =\,diag(I_l,0)$.
Now we just take $P = diag(I_l,0) A^{-1}$, then $Q = P^TP$. 
\end{proof}
\begin{Rque}
This proposition is the principal result in \cite{Mos} (theorem 2), proved by a much lengthier argument. Our proof works also for symmetric matrices
over an algebraically closed field. 
\end{Rque}

Now let $Q$ be a non-degenerate symmetric bilinear form on $V$. Recall that a $k$-plane $\Lambda \subset V$ is said {\em isotropic} if 
$Q|_{\Lambda \times \Lambda} =0$. Let $\Sigma_k$ be the locus of isotropic subspaces of $Q$. It is well-known that $\Sigma_k$
is not empty iff $k \leq m/2$, and irreducible if $k < m/2$. Let $$G = \{A \in GL(V)\, | \,Q(Av,Aw) = Q(v,w) \, \,for\, all \, v,w \in V\},$$ which acts naturally on 
$\Sigma_k$. The main purpose is to prove the following well-known fact (see \cite{ACGH}, p102):
\begin{Prop}
The action of $G$ on $\Sigma_k$ is transitive.
\end{Prop}
\begin{proof}
As easily seen, if $k = m/2$, then this proposition is true. In the following, we suppose $k < m/2$, thus $\Sigma_k$ is irreducible.
Since $Q$ is non-degenerate, our previous proposition shows that we can suppose $Q = I_m$. Now $G = \{A \in GL(V) | A^TA =I \}$.
Any $k$-plane is represented by an $m \times k$ matrix $\Lambda$. We can suppose 
$$ \Lambda = \begin{pmatrix}I_k \\  M_1 \\ M_2 \end{pmatrix}$$ where $M_1$ is a $k\times k$ matrix, and $M_2$ is an $(m-2k) \times k$ matrix.
The isotropic condition is $$I_k + M_1^TM_1 + M_2^TM_2 = 0. $$
 By this we see that $rk(M_1^T, M_2^T) = k$. Using some matrix of the form $I_m + E_{ij}+E_{ji} -E_{ii} -E_{jj} \in G$, where
$E_{ij}$ is the matrix with 1 at the position $(i,j)$ and 0 elsewhere, we can suppose $M_1$ is invertible. \\
\quad \\
{\bf Key claim:} For any $\Lambda \in \Sigma_k$, the orbit $G \cdot \Lambda$ is dense in $\Sigma_k$.

\begin{proof}[Proof of the key claim]
Take a matrix            
\begin{equation*}  A= \begin{pmatrix} I_k & 0 & 0 \\ 0 & A_1 & A_2 \\ 0 & B_1 & B_2 \end{pmatrix} \end{equation*}
we will find
 out some conditions for $\Lambda'$
such that there exists some $A \in G$ of the above form with $A \Lambda = \Lambda'$.  The last condition gives 
$$A_1 = N_1 - A_2 M , \qquad B_1 = N_2 - B_2 M $$ where $ N_1 = M'_1M_1^{-1},  N_2 = M'_2M_1^{-1},
 M = M_2 M_1^{-1}.$  The isotropic conditions for $\Lambda, \Lambda'$ imply that $$N_1^TN_1 + N_2^TN_2 = I_k + M^TM.$$

The condition $A^TA=I_m$ gives three equations. Combining with the above equation, we see that among the three there are only two independent ones, 
which can be simplified as 
\begin{equation*} \left \{ \begin{aligned}  N_1^T A_2 + N_2^T B_2 &= M^T \\ A_2^T A_2 + B_2^T B_2 &= I_{m-2k}.\end{aligned} \right. \end{equation*}

From the first equation, we can resolve for $A_2$, and the only equation for $B_2$ reads $$B_2^T X B_2 - Y^T B_2 - B_2^T Y + Z = I_{m-2k}, $$ 
where $X, Z$ are symmetric, $X = I_{m-2k} + N_2 N_1^{-1} (N_1^T)^{-1} N_2^T$, so the condition $rk(X) = m-2k$ (maximal possible rank) is an open condition on
$\Lambda' \in \Sigma_k$. If this condition is satisfied, our previous proposition gives $X = P^T P$ for some invertible matrix $P$.
Now the equation for $B_2$ is equivalent to $$ (B_2^T P^T - Y^TP^{-1})(PB_2 - (P^{-1})^TY) = W $$ where $W$ is symmetric, thus of the form $V^TV$. Now we 
just take $PB_2 - (P^{-1})^TY = V$, i.e. $B_2 = P^{-1}V + X^{-1}Y$.   So if $\Lambda'$ satisfies 
$$rk(I_{m-2k} + M'_2 M'_1{}^{-1} (M'_1{}^T)^{-1} M'_2{}^T) = m-2k,$$ then $\Lambda' \in G\cdot \Lambda$.
As $\Sigma_k$ is irreducible, we conclude that $G \cdot \Lambda$ is dense in $\Sigma_k$.   
\end{proof}

To finish the proof of the proposition, we just note that there exists at least one closed orbit, say $G \cdot \Lambda$. Then our key claim implies that 
$G \cdot \Lambda = \Sigma_k$, thus the action of $G$ on $\Sigma_k$ is transitive.
\end{proof}

{\em Acknowledgment}: This proof is done during the workshop on Algebraic Curves at UNSA, Nice. The author wants to thank J. Kock for reading a first 
draft of this note.

\quad \\
Labortoire J.A.Dieudonn\'e, Parc Valrose \\ 06108 Nice cedex 02, FRANCE \\
baohua.fu@polytechnique.org
\end{document}